\documentclass[reqno]{amsart}
\usepackage{amsmath,amsthm,amssymb,hyperref,graphicx,mathrsfs,mathtools} 

\newtheorem{theorem}{Theorem}
\newtheorem{corollary}{Corollary}

\newtheorem{thmx}{Theorem}
\newtheorem{remark}{Remark}

\usepackage[left=1in, right=1in]{geometry}

\allowdisplaybreaks

\begin{document}

\title[ON THE LOCATION OF  ZEROS OF A QUATERNIONIC POLYNOMIAL ]{ON THE LOCATION OF THE ZEROS OF A QUATERNIONIC POLYNOMIAL}
\author{N. A. Rather$^{1}$}
\author{Tanveer Bhat$^2$}
\address{$^{1,2}$ Department of Mathematics, University of Kashmir, Srinagar-190006, India}
\email{$^1$ dr.narather@gmail.com, $^2$tanveerbhat054@gmail.com}  

\begin{abstract}
In this paper, we are concerned with the problem of locating the zeros of polynomials of a quaternionic variable with quaternionic coefficients. We derive some new Cauchy bounds for the zeros of a polynomial by virtue of maximum modulus theorem.  Our results will generalise some recently proved results about the distribution of zeros of a quaternionic polynomial.
\end{abstract}

\maketitle

\section{\textbf{INTRODUCTION }}
In order to introduce the framework in which we will work, let us first introduce some preliminaries on quaternions which will be useful in the sequel. Quaternions are a unique and intriguing mathematical concept that extends the idea of complex numbers into four dimensional space. They were first introduced by the Irish mathematician Sir William Rowan Hamilton (1805-1865) in 1843. Hamilton became obsessed with quaternions and their uses \cite{aa} after inventing them and he did so for the remainder of his life. However, he probably never imagined that his invention quaternions would one day be used to programme video games and steer spacecraft. Unlike complex numbers, which have one real part and one imaginary part, quaternions have three imaginary parts. They are typically represented in the form $\mathbb{H}=\{{a=a_{0}+a_{1}i+a_{2}j+a_{3}k:a_{0},a_{1},a_{2},a_{3}\in\mathbb{R}\\}$ and $i, j, k$ satisfy $i^2=j^2=k^2=ijk=-1, ~ ij=-ji=k,~ jk=-kj=i,~ ki=-ik=j.$ One of the key properties of quaternions is non commutativity that makes them distinct from real and complex numbers. Quaternions find applications in various fields including computer graphics, robotics, quantum mechanics and theoretical physics. In computer graphics for instance they are used to represent relations in three dimensional space due to their compact representation and avoidance of singularities. Overall quaternions are a fascinating mathematical structure with wide ranging applications and connections to diverse area of mathematics and physics. Their unique properties and rich structure continue to fascinate mathematicians and scientists alike, making them an important subject of study in modern mathematics. \\
The number in the quaternions is denoted by q where $q = {\alpha}+{\beta}i+{\gamma}j+{\delta}k \in\mathbb{H}$. The conjugate of $q$, denoted by $\bar{q}$ is a quaternion $q=\alpha-\beta i-\gamma j-\delta k$ and the norm of a quaternion is defined as the square root of the sum of squares of its components, denoted by $||q||=\sqrt{q\bar{q}}=\sqrt{{\alpha}^2+{\beta}^2+{\gamma}^2+{\delta}^2}$.  The quaternion with unit norm are known as unit quaternions. The inverse of each non zero element $q$ of $\mathbb{H}$ is given by $q^{-1}=|q|^{-2}\bar{q}$.\\
Depending upon the position of the coefficients, the quaternion polynomial of degree n  in indeterminate $q$ is defined as 
$f(q) = q^{n} +q^{n-1}a_{1}+ \cdots +qa_{n-1} + a_n$ or $g(q) = q^{n} +a_{1}q^{n-1}+ \cdots +a_{n-1}q + a_n.$\\

 The exploration of the zeros a polynomial  has a deeply intricate past, serving as a wellspring of inspiration for both practical applications and theoretical inquiries, notably fuelling the development of modern algebra. While methods for determining polynomial roots using algebraic and analytic approaches can often be complicated, imposing certain constraints on polynomials proves advantageous. This field traces back to the era coinciding with the introduction of geometric representations of complex numbers in mathematics, with pioneering contributions emerging from personalities like Gauss and Cauchy. In 1829, A. L. Cauchy provided a straightforward formulation regarding the location of zeros for a polynomial with complex coefficients. He presented a concise expression defining the area encompassing all the zeros based solely on the polynomial’s coefficients. In fact he proved the following theorem:
\begin{thmx}\label{od}
If $f(z) = z^{n} +a_{n-1}z^{n-1}+...+a_1z + a_0$ is a complex polynomial, then all the zeros of $p(z)$ lie inside the disc $\left|z\right| < 1 + \max\limits_{0 \leq \nu \leq {n-1}} |a_{\nu}| $.
\end{thmx}
As a generalisation of Cauchy's classical result Kuniyeda \cite{KM}, Montel \cite{MP}, and Toya \cite{TT} proved the following result:
\begin{thmx}\label{pd}
For any $r$ and $s$, such that $r>1$, $s>1$, \quad $\frac{1}{r} + \frac{1}{s} = 1$, the polynomial $f(z) = a_{o} + a_{1}z +\cdots + z^{n}$ has all its zeros in the circle
\begin{equation}\nonumber
\begin{split}
\left|z\right| &\leq \left\lbrace 1 + \left[\sum\limits_{j=0}^{n-1} \left|a_{j}\right|^{r}\right]^{\frac{s}{r}}\right\rbrace^{\frac{1}{s}}\\[1em]
 & \leq \left(1 + n^{\frac{s}{r}} M^{s}\right)^{\frac{1}{s}} 
\end{split}
\end{equation}
Where $M = max \left|a_{j}\right|, \quad j = 0, 1, \cdots, {n-1}$.
\end{thmx}
Further as a generalisation of Cauchy's theorem Montel \cite{MR} proved the following result:
\begin{thmx}\label{qd}
All the zeros of the polynomial $f(z) = z^{n} + a_{n-1}z^{n-1} + \cdots + a_{z} + a_{0}$ lie in
$\left|z\right| \leq max\left(L, L^{\frac{1}{n}}\right), \\ \quad \text{where} \quad L = \left|a_{0}\right| + \left|a_{1}\right| + \cdots + \left|a_{n-1}\right|.$
\end{thmx}
Dar et al. \cite{RD} proved the quaternionic version of Cauchy's theorem by proving:
\begin{thmx}\label{rd}
If $f(q) = q^{n} + q^{n-1}a_{n-1} + q^{n-2}a_{n-2} +...+ qa_{1}+a_{0}$ is a quaternion polynomial with quaternion coefficients and $q$ is quaternionic variable, then all the zeros of $f(q)$ lie in the ball $\left|q\right| < 1 + \max\limits_{0 \leq \nu \leq {n-1}} \left|a_{\nu}\right|$.  
\end{thmx}
Further various refinements and generalizations of Cauchy's theorem in quaternion settings were proved in subsequent years. In this direction Rather et al. \cite{NR} proved various results regarding the location of zeros of quaternionic polynomials with quaternionic coefficients, besides they refined the cauchy's theorem in quaternionic settings by proving the following result:
\begin{thmx}\label{sd}
Let $f(q) = q^{n} + q^{n-1}a_{1} +\cdots + qa_{n-1} + a_{n}$ be a quaternion polynomial of degree $n$ with quaternion coefficients and $q$ be a quaternion variable. If $\alpha_{2} \geq \alpha_{3} \geq \cdots \geq \alpha_{n}$ are ordered positive numbers,
\begin{equation}
\alpha_{\nu} = \frac{\left|a_{\nu}\right|}{r^{\nu}}, \quad \nu = 2, 3, \cdots,n
\end{equation}
where $r$ is a positive real number. Then all the zeros of $f(q)$ lie in the union of balls\\
$\left\lbrace q\in H : \left|q\right| \leq r\left(1 + \alpha_{2}\right)\right\rbrace \quad and \quad \left\lbrace q\in H : \left|q+ a_{1}\right| \leq r\right\rbrace$.
\end{thmx} 
\section{\textbf{ Main results}}
Unlike polynomials over real and complex numbers, the zeros of quaternion polynomials can exhibit more intricate behaviour due to the non commutativity of quaternion multiplication. Research in this area often focusses on developing methods for locating and characterizing the zeros of quaternionic polynomials. This may involve extending techniques from complex analysis or algebra to quaternionic settings as well as developing new tools specific to quaternionic polynomials. Now regarding the location of zeros of a quaternion polynomial with quaternion coefficients, we begin by proving the quaternionic version of Theorem \ref{pd} and Theorem \ref{qd}.
\begin{theorem}\label{L2}
For any $r$ and $s$ such that $r>1$, $s>1$, $\frac{1}{r} + \frac{1}{s} = 1$, let $f(q) = q^{n} + a_{n-1}q^{n-1} + \cdots + a_{1}q + a_{0}$ be a quaternion polynomial with quaternionic coefficients and q be quaternionic variable, then all the zeros of $f(q)$ lie in the ball
\begin{align}\label{Ee2}
\left|q\right| &\leq \left\lbrace1 + \left[\sum\limits_{\nu=0}^{n-1} |a_{\nu}|^{r}\right]^{\frac{s}{r}}\right\rbrace^{\frac{1}{s}}\\[1em]
\begin{split}\label{Ee3}
\left|q\right| &\leq \left[1 + n^{\frac{s}{r}}M^{s}\right]^{\frac{1}{s}}.
\end{split}
\end{align} 
\end{theorem} 
Where $M  = max\left|a_{\nu}\right|, \quad \nu = 0,1, \cdots, {n-1}.$
\begin{remark}
If we let $r \to \infty$ so that $s \to 1$, Theorem \ref{L2} reduces to the quaternion version of Cauchy's theorem due to Dar et al. $\left(\text{Theorem} \ \ref{rd} \right)$.
\end{remark}
\begin{corollary}\label{c1}
If we take $r = s = 2$ in Theorem \ref{L2}, then all the zeros of $f(q)$ lie in the ball
\begin{equation}
 \left|q\right| \leq \left\lbrace 1 + \sum\limits_{\nu=0}^{n-1} |a_{\nu}|^{2}\right\rbrace^{\frac{1}{2}}.
\end{equation}
\end{corollary}
\begin{theorem}\label{tA}
Let $f(q) = q^{n} + a_{n-1}q^{n-1} + \cdots + a_{1}q + a_{0}$  be a quaternion polynomial of degree n with quaternion coefficients, then all the zeros of $f(q)$ lie in the ball\\
\begin{equation}
\left|q\right| \leq max \left( L , L^{\frac{1}{n}}\right),\quad  \text{where} \quad L = \left|a_{0}\right| + \left|a_{1}\right| +\cdots + \left|a_{n-1}\right|. 
\end{equation}
\end{theorem}
\section{\textbf{Proof of the main theorems}}
\begin{proof}[\bf{Proof of Theorem \ref{L2}}]
Let $f(q) = q^{n} + a_{n-1}q^{n-1} + \cdots + a_{1}q + a_{0}$, therefore\\
\begin{align*}
\left|f(q)\right| \geq \left|q\right|^{n}\left[ 1 - \left\lbrace \left|a_{n-1}\right|\frac{1}{\left|q\right|} + \cdots + \left|a_{0}\right|\frac{1}{\left|q\right|^{n}}\right\rbrace\right]
\end{align*}\
\begin{align*}
 = \left|q\right|^{n}\left\lbrace 1 - \sum\limits_{\nu=0}^{n-1} \left|a_{\nu}\right|\frac{1}{\left|q\right|^{n-\nu}}\right\rbrace. 
\end{align*}
Using Holder's inequality, we get
\begin{align*}
 \left|f(q)\right| \geq \left|q\right|^{n}\left[ 1 -  \left\lbrace\sum\limits_{\nu=0}^{n-1} \left|a_{\nu}\right|^{r}\right\rbrace^{\frac{1}{r}}\left\lbrace\sum\limits_{\nu=0}^{n-1} \frac{1}{\left|q\right|^{\left(n-\nu\right)s}}\right\rbrace^{\frac{1}{s}}\right] 
\end{align*}
\begin{align}\label{t1e1}
&= \left|q\right|^{n}\left[ 1 -  \left\lbrace\sum\limits_{\nu=0}^{n-1} \left|a_{\nu}\right|^{r}\right\rbrace^{\frac{1}{r}}\left\lbrace\sum\limits_{\nu=0}^{n} \frac{1}{\left|q\right|^{\nu s}}\right\rbrace^{\frac{1}{s}}\right].
\end{align}
Let $A_{r} = \left\lbrace\sum\limits_{\nu=0}^{n-1}\left|a_{\nu}\right|^{r}\right\rbrace^{\frac{1}{r}}$. If $\left|q\right|>1$, we have
\begin{equation}\nonumber
\begin{split}
\sum\limits_{\nu=0}^{n} \frac{1}{\left|q\right|^{\nu s}} &< \sum\limits_{\nu=0}^{\infty} \frac{1}{\left|q\right|^{s}} \\
&= \frac{\frac{1}{\left|q\right|^{s}}}{1 - \frac{1}{\left|q\right|^{s}}}\\
&= \frac{1}{\left|q\right|^{s} - 1}.
\end{split}
\end{equation}
Using this in (\ref{t1e1}), we get for $\left|q\right| > 1$
\begin{equation}\nonumber
\begin{split}
\left|f(q)\right| &> \left|q\right|^{n}\left\lbrace 1 - \frac{A_{r}}{\left(\left|q\right|^{s}-1\right)^{\frac{1}{s}}}\right\rbrace \\
&> 0 \quad \text{if}  \quad \left|q\right|^{s}-1 \geq \left(A_{r}\right)^{s}.\\
\text{i,e if} \quad \left|q\right| &\geq \left( 1 + \left(A_{r}\right)^{s}\right)^{\frac{1}{s}}.
\end{split}
\end{equation}
This shows that those zeros of $f(q)$ which lie outside the unit ball lie in \ $\left|q\right| \leq \left[1 + \left(A_{r}\right)^{s}\right]^{\frac{1}{s}}$.
But those zeros of $f(q)$, which lie inside or on the unit ball already satisfy inequality \ $\left|q\right| \leq \left[1 + \left(A_{r}\right)^{s}\right]^{\frac{1}{s}}$.
Therefore all the zeros of $f(q)$ lie in
\begin{align*}
\left|q\right| \leq \left\lbrace1 + \left[\sum\limits_{\nu=0}^{n-1} |a_{\nu}|^{r}\right]^{\frac{s}{r}}\right\rbrace^{\frac{1}{s}}.
\end{align*}
Now $ M = \max\limits_{0\leq \nu \leq {n-1}}\left|a_{\nu}\right|$,\
this gives \quad $\left|a_{\nu}\right|^{r} \leq M^{r},\quad  \nu = 0,1, \cdots, {n-1}$. \\
Implies
\begin{align*}
\sum\limits_{\nu=0}^{n-1}\left|a_{\nu}\right|^{r} \leq nM^{r}.
\end{align*}
Therefore
\begin{align*}
\left\lbrace\sum\limits_{\nu=0}^{n-1} \left|a_{\nu}\right|^{r}\right\rbrace^{\frac{s}{r}} \leq n^{\frac{s}{r}}M^{s}.
\end{align*}
And hence from (\ref{Ee2}), we get all the zeros of $f(q)$ lie in \ 
$\left|q\right| \leq \left[1 + n^{\frac{s}{r}}M^{s}\right]^{\frac{1}{s}}$,\ 
Which is (\ref{Ee3}).\\
This completes the proof of Theorem \ref{L2}.
\end{proof}
\begin{proof}[\bf{ Proof of Theorem \ref{tA}}]
We have $f(q) = q^{n} + a_{n-1}q^{n-1} + \cdots + a_{1}q + a_{0}$, then
\begin{align*}
\left|f(q)\right| \geq \left|q\right|^{n} - \left|a_{n-1}q^{n-1} + \cdots + a_{1}q + a_{0}\right|
\end{align*}
\begin{align}\label{Ee7}
\geq \left|q\right|^{n}\left\lbrace 1 - \sum\limits_{\nu=0}^{n}\left|a_{n-\nu}\right|\frac{1}{\left|a\right|^{\nu}}\right\rbrace.
\end{align}
Let $L\geq1$, then $max\left(L, L^{\frac{1}{n}}\right) = L$.\\
Therefore for $\left|q\right| \geq 1$, we have \quad 
$\left|q\right|^{\nu} \geq \left|q\right|, \quad \nu = 1,2, \cdots, n$.
Hence from equation (\ref{Ee7}), we have
\begin{equation}\nonumber
\begin{split}
\left|f(q)\right| &\geq \left|q\right|^{n} \left\lbrace 1 - \frac{1}{\left|q\right|}\sum\limits_{\nu=1}^{n}\left|a_{n-\nu}\right|\right\rbrace\\
& = \left|q\right|^{n}\left\lbrace 1 - \frac{1}{\left|q\right|}L\right\rbrace\\
&> 0
\end{split}
\end{equation}
if $1 - \frac{L}{\left|q\right|} > 0$,\quad
that is  if $\left|q\right| > L$.
That is $f(q)$ does not vanish for $\left|q\right| > L$.
Hence all the zeros of $f(q)$ which lie outside of the unit ball lie in $\left|q\right| \leq L = max \left(L, L^{\frac{1}{n}}\right)$ in case $L \geq 1$. But those zeros of $f(q)$ which lie inside the unit ball already lie in $\left|q\right| \leq L$.
Therefore all the zeros of $f(q)$ lie in $\left|q\right| \leq L = max\left(L, L^{\frac{1}{n}}\right)$ in case $L \geq 1$.\\
Now suppose that $L < 1$, then $max\left(L, L^{\frac{1}{n}}\right) = L^{\frac{1}{n}}$, and for $\left|q\right| < 1$, we have $\left|q\right|^{n} < \left|q\right|^{\nu},\quad \nu = 1,2,\cdots,n$.\\
So that $-\frac{1}{\left|q\right|^{n}} \leq - \frac{1}{\left|q\right|^{\nu}}, \quad \nu = 1,2, \cdots,n$.
Using this in (\ref{Ee7}), we get
\begin{equation}\nonumber
\begin{split}
\left|f(q)\right| &\geq \left|q^{n}\right| \left\lbrace1 - \frac{1}{\left|q\right|^{n}}\sum\limits_{\nu=1}^{n} \left|a_{n-\nu}\right|\right\rbrace \\
&= \left|q\right|^{n}\left\lbrace1 - \frac{1}{\left|q\right|^{n}}L\right\rbrace\\
&> 0
\end{split}
\end{equation}
if $\left|q\right| > L$  or if $\left|q\right| > L^{\frac{1}{n}}$.
That is $f(q)$ does not vanish for $\left|q\right| > L^{\frac{1}{n}}$.
Therefore all the zeros of $f(q)$ lie in $\left|q\right| \leq L^{\frac{1}{n}} = max\left( L, L^{\frac{1}{n}}\right)$ in case $L <1$. Hence it follows that all the zeros of $f(q)$ lie in $\left|q\right| \leq max\left( L, L^{\frac{1}{n}}\right)$, Where $L = \left|a_{0}\right| + \left|a_{1}\right| + \cdots + \left|a_{n-1}\right| = \sum\limits_{\nu=1}^{n}\left|a_{n-1}\right|$.\\
This completes the proof of Theorem \ref{tA}.
\end{proof}
\section{Declarations}
\subsection{Availability of data and material}
Not applicable.
\subsection{Competing interests}
The authors declare that they have no competing interests.
\subsection{Funding}
None.
\subsection{Author's contributions}
All authors contributed equally to the writing of this paper. All authors read and approved the final manuscript.

\end{document}